\documentclass[11pt]{article}
\usepackage[utf8]{inputenc}
\usepackage{amsmath, amssymb, amsthm}
\usepackage{geometry}
\usepackage{hyperref}
\hypersetup{
	colorlinks=true,
	linkcolor=blue,
	filecolor=magenta,      
	urlcolor=cyan,
}

\newtheorem{conjecture}{Conjecture}

\title{On the Representation of Integers as Sums of Limited Prime Powers}
\author{Julius Stricker} 
\date{}

\begin{document}
	
	\maketitle
	
	\begin{abstract}
		We present a novel conjecture concerning the additive representation of natural numbers using prime powers. Based on extensive computational verification, we conjecture that every integer $n > 23$ can be expressed as a sum of at most five prime powers $p^k$, where $p$ is a prime number and $k \ge 2$ is an integer. This conjecture is supported by comprehensive computational evidence covering all integers up to $10^7$ (exhaustively) and specific large numbers up to $10^{10}$ (via sampling), where no counterexample requiring more than five summands has been found. This work highlights a surprising "combinatorial creativity" of prime powers, which allows for efficient additive representations despite their asymptotic sparsity and the existence of extremely large gaps between individual prime powers.
	\end{abstract}
	
	\section{Introduction}
	Additive number theory is a rich field concerned with representing integers as sums of elements from specific sets. Classical problems in this area include Waring's Problem \cite{Waring}, which states that every natural number can be expressed as a sum of a fixed number of $k_0$-th powers of natural numbers, and Goldbach's Conjectures \cite{Goldbach}, which deal with sums of prime numbers.
	
	Our work explores a variation of these problems by restricting the summands to a specific, sparser set: \textit{prime powers} $p^k$ where $p$ is a prime and the exponent $k$ is an integer greater than or equal to $2$. This set includes numbers like $2^2=4, 3^2=9, 2^3=8, 5^2=25, 3^3=27$, etc. Despite the relative sparsity of this set compared to the set of prime numbers, our extensive computational analysis reveals a surprisingly efficient additive representation for natural numbers.
	
	We formally state our primary finding as a conjecture for all numbers greater than 23.
	
	\begin{conjecture}[Main Conjecture on Limited Prime Power Sums]
		Every natural number $n > 23$ can be expressed as a sum of $m$ prime powers $P_i^{k_i}$, where $P_i$ is a prime number and $k_i \ge 2$ is an integer for each $i$, and the number of summands $m$ satisfies $2 \le m \le 5$.
	\end{conjecture}
	
	This conjecture is particularly intriguing due to its counter-intuitive nature when compared to established results and conjectures concerning sums of prime numbers.
	
	\section{Definitions}
	A \textbf{prime power} is a positive integer of the form $p^k$, where $p$ is a prime number and $k$ is an integer such that $k \ge 2$.
	Examples of prime powers include: $2^2=4, 3^2=9, 2^3=8, 5^2=25, 3^3=27, 2^4=16, 7^2=49, \dots$.
	Note that a prime number $p$ itself (i.e., $p^1$) is not considered a prime power under this definition.
	
	\section{Computational Evidence and Methodology}
	The foundation of this work lies in extensive computational verification. Our methodology involved systematically attempting to represent integers as sums of prime powers.
	
	\subsection{Support for Conjecture 1}
	For Conjecture 1, which extends to all numbers $n > 23$ with at most five summands, we have gathered substantial numerical evidence:
	\begin{itemize}
		\item All integers up to $10^7$ have been exhaustively checked. This process involved an algorithm that, for each integer $n$ in this range, searches for a combination of up to five prime powers that sum to $n$. No counterexample requiring more than five summands was found up to $10^7$.
		\item Further computational checks were performed on $100,000$ numbers starting from $1,500,000$ and extending up to $10^{10}$. These large numbers were also successfully represented as sums of at most five summands.
\item Analysis of prime power distribution reveals that the gaps between consecutive prime powers can be extremely large. For instance, the gap between the squares of two consecutive prime numbers, $p_n = 18,361,375,334,787,046,697$ and $p_{n+1} = 18,361,375,334,787,048,269$, amounts to over $57.7$ Quintillion ($p_{n+1}^2 - p_n^2 = 57,728,164,052,570,477,286,552$). Despite these immense gaps between individual prime powers, our computational evidence strongly suggests that a small number of summands is sufficient.
	\end{itemize}
	Due to resource limitations, exhaustive checks for numbers beyond $10^7$ were not feasible, and larger ranges were explored via targeted sampling. The consistent success across these vast ranges provides strong empirical support for the conjecture.
	
	\subsection{Examples of Representations}
	The following examples illustrate the representations, showcasing the variety of prime bases and exponents used:
	\begin{itemize}
		\item $24 = 2^3 + 2^4 = 8 + 16$ (2 summands)
		\item $26 = 3^2 + 3^2 + 2^2 + 2^2 = 9 + 9 + 4 + 4$ (4 summands)
		\item $78869 = 263^2 + 19^3 + 53^2 + 2^5 = 69169 + 6859 + 2809 + 32$ (4 summands)
		\item $512550 = 502681 + 7921 + 1331 + 361 + 256$ (5 summands)
		\item $701721 = 687241 + 12769 + 961 + 625 + 125$ (5 summands)
		\item $2500050 = 2493241 + 5329 + 625 + 512 + 343$ (5 summands)
		\item $12000728 = 11992369 + 6889 + 841 + 625 + 4$ (5 summands)
		\item $50000162 = 49970761 + 24389 + 2809 + 2187 + 16$ (5 summands)
		\item $250000558 = 249924481 + 73441 + 1369 + 1024 + 243$ (5 summands)
		\item $2500000710 = 2499900001 + 96721 + 2401 + 1331 + 256$ (5 summands)		
		\item $3000000001 = 2999424289 + 552049 + 15625 + 4913 + 3125$ (5 summands)
		
	\end{itemize}
	
	\section{Discussion and Significance}
	The most striking aspect of this conjecture is its counter-intuitive nature when compared to established results in additive number theory.
	
	\subsection{Comparison with Goldbach's Conjectures}
	The Strong Goldbach Conjecture \cite{Goldbach} states that every even number greater than 2 is the sum of two primes. The Weak Goldbach Conjecture \cite{Vinogradov} (now a theorem by Helfgott \cite{Helfgott}) states that every odd number greater than 5 is the sum of three primes. Prime numbers are asymptotically denser than prime powers $p^k$ with $k \ge 2$.
	
	Intuitively, one might expect that a sparser set of building blocks would require more summands to cover all numbers. However, our findings present a fascinating contrast:
	\begin{itemize}
		\item For prime numbers (a denser set), we need up to 3 summands to cover all numbers (if Strong Goldbach holds).
		\item For prime powers $p^k, k \ge 2$ (a much sparser set), our computational evidence suggests that 5 summands suffice for all numbers $n > 23$.
	\end{itemize}
	This behavior, while initially counter-intuitive, becomes logically explainable when considering the unique properties of prime powers. While prime numbers are "stiff" (always $p^1$), prime powers offer "combinatorial flexibility" through varying exponents ($k=2, 3, 4, \dots$). This allows for an exponential increase in the number of possible combinations as the target number grows, effectively compensating for the lower asymptotic density of individual prime powers. It is as if the race starts in "rear gear" compared to primes, but then accelerates to "light speed" towards infinity due to the sheer volume of available sum configurations. This suggests that the "combinatorial creativity" of prime powers can bridge gaps of any kind, leading to surprisingly efficient representations.
	
	\subsection{The Density of Primes vs. Prime Powers}
	By the prime number theorem (PNT), the number of primes $\pi(x)$ up to $x$ is approximately $\frac{x}{\ln(x)}$. This principle can be extended to analyze the density of prime powers. The number of $k$-th powers of a prime below $x$ is of the order $x^{1/k}/\ln(x^{1/k})$.
	
	Specifically, the number of prime squares ($p^2 \le x$) is approximately $\pi(\sqrt{x}) \approx \frac{\sqrt{x}}{\ln(\sqrt{x})}$, while the number of prime cubes ($p^3 \le x$) is approximately $\pi(\sqrt[3]{x}) \approx \frac{\sqrt[3]{x}}{\ln(\sqrt[3]{x})}$. The total number of prime powers $p^k$ with $k \ge 2$ up to $x$ is therefore given by the sum $\sum_{k=2}^{\lfloor \log_2 x \rfloor} \pi(x^{1/k})$. The term for prime squares dominates this sum, as the other terms grow significantly slower with $x$.
	
	From this density, the number of $t$-term expressions with sums $\le x$ is a combinatorial product of these densities. For example, the number of ways to form a sum of two prime squares below $x$ is approximately $\frac{1}{2!} \left( \frac{\sqrt{x}}{\ln(\sqrt{x})} \right)^2$. While this provides an intuitive measure of the search space, it does not account for the 'combinatorial creativity' of using different exponents ($p_1^{k_1} + p_2^{k_2} + \dots$). Our computational evidence demonstrates that this flexibility is key to covering the entire number line with a small number of summands.
	
	\subsubsection{Comparison of Densities}
	To illustrate the stark difference in density, the following table compares the exact number of prime numbers $\pi(n)$ with the number of all prime powers $p^k \le n$ (where $k \ge 2$) for various values of $n$.
	
	\begin{center}
		\begin{tabular}{|c|c|c|}
			\hline
			\textbf{Upper Limit $n$} & \textbf{Number of Primes $\pi(n)$} & \textbf{Number of Prime Powers ($p^k \le n, k \ge 2$)} \\
			\hline
			$10^4$ & $1229$ & 51 \\
			\hline
			$10^6$ & $78498$ & 236 \\
			\hline
			$10^8$ & $5761455$ & 1404 \\
			\hline
		\end{tabular}
	\end{center}
	
	\subsection{Relationship to the Waring-Goldbach Problem}
	Our conjecture can be seen as a unique variant of the classical Waring-Goldbach Problem, which concerns the representation of integers as sums of fixed exponent $k_0$-th powers of prime numbers.
	
	In the Waring-Goldbach context, established theorems (for "sufficiently large numbers") include:
	\begin{itemize}
		\item For $k_0=2$ (sums of squares of primes, $p^2$): Every sufficiently large **odd** number can be expressed as a sum of at most **5** squares of primes.
		\item For $k_0=3$ (sums of cubes of primes, $p^3$): Every sufficiently large **natural number** (both even and odd) can be expressed as a sum of at most **9** cubes of primes.
	\end{itemize}
	
	Our conjecture distinguishes itself by allowing **mixed prime powers** ($p^k$ where $k \ge 2$ can vary for each summand) to represent **all natural numbers** (even and odd) greater than 23. This is a crucial difference:
	\begin{itemize}
		\item While the $p^2$ and $p^3$ terms are numerically the most dominant and frequently used prime powers in our sums, our set of allowed summands is the **union of all $k$-th powers of primes for $k \ge 2$**. This includes $p^4, p^5$, etc., offering a significantly richer set of "building blocks" compared to the fixed-exponent approach of the classical Waring-Goldbach problem.
		\item This flexibility, allowing the exponent $k$ to vary, introduces a **higher degree of combinatorial complexity** for analytical proofs. Traditional methods, which often rely on the uniform structure of fixed exponents, become substantially more challenging when dealing with such a mixed set. This increased complexity is likely a reason why this specific problem (sum of mixed prime powers) has not been as prominently studied as its fixed-exponent counterparts.
		\item Despite this added complexity in the set of summands, our computational evidence suggests a remarkably low upper bound of 5 summands. This highlights the immense "combinatorial creativity" of these mixed prime powers, enabling them to cover the entire number line efficiently, even across the vast gaps between individual prime powers.
	\end{itemize}
	Thus, our conjecture represents a unique and challenging extension within the additive number theory landscape, bridging the gap between specific fixed-power problems and a more general, flexible representation.
	
	\subsection{Implications for Proof Strategy and Computational Limits}
	The extensive computational evidence provides strong empirical support for Conjecture 1. However, the analysis of prime power distribution reveals that the gaps between consecutive prime powers can be substantial, reaching trillions or even quintillions in higher ranges. 
	
	The ultimate proof of Conjecture 1 would likely require advanced methods from analytic number theory, such as the Hardy-Littlewood \cite{HardyLittlewood} circle method. These methods aim to prove the existence of representations for "sufficiently large numbers" by showing that the number of ways to represent $N$ as a sum of $m$ prime powers is non-zero for $N \to \infty$. This approach does not rely on exhaustive computational verification for all numbers, which becomes impossible given the astronomical number of combinations in higher ranges. The current computational evidence serves as strong empirical support for the conjecture, bridging the gap between small numbers and the theoretical realm of "sufficiently large numbers."

	\section{Conclusion and Future Work}
	Based on extensive computational evidence, we conjecture that every integer $n > 23$ can be expressed as a sum of at most five prime powers ($p^k, k \ge 2$). This finding challenges intuitive notions about the density of additive bases and highlights the unique combinatorial power of prime powers.
	
	Future work includes:
	\begin{itemize}
		\item A formal mathematical proof of Conjecture 1, likely employing methods from analytic number theory to establish the existence of representations for all sufficiently large numbers.
		\item Further computational verification for even larger numbers, potentially pushing the boundaries of the current evidence, although exhaustive checks for all combinations are computationally infeasible.
	\end{itemize}
	
	\section{Data Availability}
	The complete computational data supporting the current work, including the verification up to $1,000,000$ and beyond, is publicly available in a dedicated GitHub repository at:
	\url{https://github.com/PrimePowers/PrimePowers}
	This repository also contains additional sample data for larger numbers as discussed in Section 3.2.

\end{document}